\documentstyle[12pt]{article}

\begin{document}

\begin{center}
{\huge {\bf Chains of Frobenius subalgebras of $so(M)$\\[3mm] and the
corresponding twists}}\\[5mm]

{\sc Vladimir D. Lyakhovsky 
\footnote{E-mail address: lyakhovs@snoopy.phys.spbu.ru ; RFBR grant 
N 00-01-00500.}\\ Theoretical
Department, St. Petersburg State University,\\ 198904, St. Petersburg,
Russia \\ \vskip0.5cm }

{\sc Alexander Stolin \footnote{E-mail address: astolin@math.chalmers.se}\\
Department of Mathematics, University of Goteborg,\\ S-412 96 Goteborg,
Sweden}\\ \vskip0.5cm

{\sc Petr P. Kulish
\footnote{E-mail address: kulish@pdmi.ras.ru \,;
On leave of absence from Steklov Mathematical Institute, Fontanka 27,
191011, St.Petersburg, Russia ; RFBR grant N 99-01-00101 
and INTAS grant N 99-01459.}\\ U. C. E. H., University of Algarve, 8000
Faro, Portugal} \\
\end{center}

\vskip0.5cm

\begin{abstract}
Chains of extended jordanian twists are studied for the universal enveloping
algebras $U(so(M))$. The carrier subalgebra of a canonical chain ${\cal F}_{%
{\cal B}_{0\prec p\max }}$ cannot cover the maximal nilpotent subalgebra $%
N^+(so(M))$. We demonstrate that there exist other types of Frobenius
subalgebras in $so(M)$ that can be large enough to include $N^+(so(M))$. The
problem is that the canonical chains ${\cal F}_{{\cal B}_{0\prec p}}$ do not
preserve the primitivity on these new carrier spaces. We show that this
difficulty can be overcome and the primitivity can be restored if one
changes the basis and passes to the deformed carrier spaces. Finally the
twisting elements for the new Frobenius subalgebras are explicitly
constructed. This gives rise to a new family of universal $R$-matrices
for orthogonal algebras. For a special case of $g = so(5)$ and its defining
representation we present the corresponding matrix solution of the
Yang-Baxter equation.
\end{abstract}

\section{Introduction}

Quantizations of triangular Lie bialgebras ${\bf L}$ with antisymmetric
classical $r$-matrices $r=-r_{21}$ satisfying the classical Yang-Baxter
equation (CYBE) form an important class of triangular Hopf algebras ${\cal A}
(m,\Delta ,S,\eta ,\epsilon ;{\cal R})$, with $R$-matrix satisfying the
unitarity condition ${\cal R}_{21}{\cal R}=1$. These quantizations are
defined by the twisting elements ${\cal F}=\sum f_{\left( 1\right) }\otimes
f_{\left( 2\right) }\in {\cal A}\otimes {\cal A}$ that satisfy the twist
equations \cite{DR83}: 
\begin{equation}
\label{twist-eq}\left( {\cal F}\right) _{12}\left( \Delta \otimes {\rm id}
\right) {\cal F}=\left( {\cal F}\right) _{23}\left( {\rm id}\otimes \Delta
\right) {\cal F},\qquad \left( \epsilon \otimes {\rm id}\right) {\cal F}
=\left( {\rm id}\otimes \epsilon \right) {\cal F}=1. 
\end{equation}

The knowledge of the twisting element is quite important in applications
giving (twisted) $R$-matrix ${\cal R_F}={\cal F}_{21}{\cal R}{\cal F}
^{-1}$ and twisted coproduct $\Delta _{{\cal F}}={\cal F}\Delta {\cal F}
^{-1} $.

The explicit expressions of the twisting elements ${\cal F}$ were found in 
\cite{KLM}, for the carrier algebras ${\bf L}$ with special properties of
their triangular decompositions. Such carrier subalgebras are the
multidimensional analogs of the enlarged Heisenberg algebra and can be found
in any simple Lie algebra $g$ of rank greater than 1. In the root system $%
\Lambda \left( g\right) $ one can choose the initial root $\lambda _0$ and
consider the set $\pi $ of its constituent roots%
$$
\begin{array}{l}
\pi =\left\{ \lambda ^{\prime },\lambda ^{\prime \prime }\mid \lambda
^{\prime }+\lambda ^{\prime \prime }=\lambda _0;\quad \lambda ^{\prime
}+\lambda _0,\lambda ^{\prime \prime }+\lambda _0\notin \Lambda \left(
g\right) \right\}  \\ 
\pi =\pi ^{\prime }\cup \pi ^{\prime \prime };\qquad \pi ^{\prime }=\left\{
\lambda ^{\prime }\right\} ,\pi ^{\prime \prime }=\left\{ \lambda ^{\prime
\prime }\right\} .
\end{array}
$$
The subalgebra ${\bf L}$ is generated by the elements $\left\{ E_{\lambda
_0},E_\lambda \mid \lambda \in \pi \right\} $ and the Cartan generator $%
H_{\lambda _0}$ dual to $\lambda _0$. The solution ${\cal F_{EJ}}$ of the
twist equations corresponding to the carrier subalgebra ${\bf L}$ is called
the extended jordanian twist and can be decomposed into the product of two
factors, the jordanian twist $\Phi _{{\cal J}}$ and the extension $\Phi _{%
{\cal E}}$: 
\begin{equation}
\label{ext}{\cal F_{EJ}}=\Phi _{{\cal E}}\cdot \Phi _{{\cal J}%
}=\prod_{\lambda ^{\prime }\in \pi ^{\prime }}\exp \left\{ E_{\lambda
^{\prime }}\otimes E_{\lambda _0-\lambda ^{\prime }}e^{-\frac 12\sigma
_{\lambda _0}}\right\} \cdot \exp \{H_{\lambda _0}\otimes \sigma _{\lambda
_0}\}.
\end{equation}
Here $\sigma _{\lambda _0}=\ln (1+E_{\lambda _0})$.

To construct the twists for the higher dimensional carrier subalgebras ${\bf %
L}$ of $g$ we have to consider the subset $\Lambda _{\lambda _0}^{\perp }$
of roots orthogonal to $\lambda _0$ and the corresponding subalgebra $%
g_{\lambda _0}^{\perp }\subset g$. It was shown in \cite{KLO} that for the
classical Lie algebras $g$ one can always find a subalgebra $g_1$ in $%
g_{\lambda _0}^{\perp }$ whose generators become primitive after the twist $%
{\cal F_{EJ}}$. Such primitivization of $g_k\subset g_{k-1}$ (called the
matreshka effect) makes it possible to compose chains of extended twists of
the type (\ref{ext}) corresponding to the injections $g_p\subset \ldots
\subset g_1\subset g_0 = g $ ,
\begin{equation}
\label{chain-ini}
\begin{array}{l}
{\cal F}_{{\cal B}_{0\prec p}}=\prod_{\lambda ^{\prime }\in \pi _p^{\prime
}}\exp \left\{ E_{\lambda ^{\prime }}\otimes E_{\lambda _0^p-\lambda
^{\prime }}e^{-\frac 12\sigma _{\lambda _0^p}}\right\} \cdot \exp
\{H_{\lambda _0^p}\otimes \sigma _{\lambda _0^p}\}\,\cdot  \\ \prod_{\lambda
^{\prime }\in \pi _{p-1}^{\prime }}\exp \left\{ E_{\lambda ^{\prime
}}\otimes E_{\lambda _0^{p-1}-\lambda ^{\prime }}e^{-\frac 12\sigma
_{\lambda _0^{p-1}}}\right\} \cdot \exp \{H_{\lambda _0^{p-1}}\otimes \sigma
_{\lambda _0^{p-1}}\}\,\,\cdot  \\ 
\ldots  \\ 
\prod_{\lambda ^{\prime }\in \pi _0^{\prime }}\exp \left\{ E_{\lambda
^{\prime }}\otimes E_{\lambda _0^0-\lambda ^{\prime }}e^{-\frac 12\sigma
_{\lambda _0^0}}\right\} \cdot \exp \{H_{\lambda _0^0}\otimes \sigma
_{\lambda _0^0}\}\,.
\end{array}
\end{equation}
In the case $g=sl(n)$ the subalgebras $g_k$ (they remain primitive after
the twisting by ${\cal F}_{{\cal E}_{k-1}{\cal J}_{k-1}}$) coincide with $%
g_{\lambda _{k-1}}^{\perp }$. The result is that the maximal canonical chain 
${\cal F}_{{\cal B}_{0\prec p\max }}$ for $g=sl(n)$ is full, its carrier
subalgebra contains all the generators of the nilpotent subalgebra $%
N^{+}\left( g\right) $.

For the orthogonal simple algebras the situation is different. In this case $%
g_{\lambda _{k-1}}^{\perp }=$ $g_k\oplus sl^{\left( k\right) }(2)$ and the
coproducts of generators in the space $g_{\lambda _{k-1}}^{\perp }\backslash
g_k$ are nontrivially deformed by the twist ${\cal F}_{{\cal E}_{k-1}{\cal J}%
_{k-1}}$. The next extended twist ${\cal F}_{{\cal E}_k{\cal J}_k}$ does not
contain these generators in its carrier space. Such chain cannot be full .

The canonical twists (\ref{chain-ini}) correspond to solutions of CYBE 
related with Frobenius subalgebras
in $g$ described by the coboundary bilinear forms $\omega
_p^{+}=\sum_{k=0}^pE_{\lambda _0^k}^{*}\left( \left[ \, , \, \right]
\right) $ \cite{STO}. 
In this paper we show that for the orthogonal algebras these forms
can be modified. The Frobenius subalgebras can be enlarged in order to
include all nonzero root generators from $g_{\lambda _{k-1}}^{\perp }
\backslash g_k$.

The problem is how to find the corresponding twists, i.e. to solve the
equations (\ref{twist-eq}) for the subalgebra $g_{\lambda _{k-1}}^{\perp }$
that contains the generators with deformed coproducts $\Delta _{{\cal F}_{%
{\cal E}_k{\cal J}_k}}$.

In \cite{KL} it was demonstrated that under certain conditions (while the
coproducts in $g$ are non-trivially twisted by ${\cal F}$ ) one can find in $%
U_{{\cal F}}(g)$ the deformed carrier subspace that is primitive and
generates a subalgebra of $g$. Below we show that this effect is in some
sense universal. The corresponding deformed spaces for orthogonal algebras
can be found for any extended twist ${\cal F}_{{\cal E}_{k-1}{\cal J}%
_{k-1}}. $As a result the canonical chain of twists ${\cal F}_{{\cal B}%
_{0\prec p}}$ can be extended using some additional factors (the deformed
jordanian twists). For the maximal value of $p$ the corresponding chain $%
{\cal F}_{{\cal G}_{0\prec p}}$ becomes full and the corresponding carrier
space contains all the generators of $N^{+}(g) $.

\section{Frobenius subalgebras}

We consider the orthogonal algebras $g=so\left( M\right) $ of the series $%
B_N $ (for $M=2N+1$) and $D_N$ (for $M=2N$). The root system $\Lambda (g)$
will be fixed as follows
\begin{equation}
\begin{array}{c}
\Lambda (g) = \left\{
\begin{array}{ll}
\pm e_i, \, \,
\pm e_j \pm e_k & {\rm for} \quad g= so(2N+1), \\ \pm e_i \pm e_j &
{\rm for} \quad g=so(2N)
\end{array}
\right\} \\
i,j,k = 1,\dots,N.
\end{array}
\end{equation}
Let ${\cal E}_{i,j}$ be the ordinary matrix units,
$$
\left[ {\cal E}_{i,j},{\cal E}_{k,l} \right]= \delta_{jk}{\cal E}_{i,l} -
\delta_{il}{\cal E}_{k,j},
$$
and $M_{i,j}$ -- the antisymmetric ones,
$$
\left[ M_{a,b},M_{c,d}\right] =\delta _{bc}M_{a,d}+\delta
_{ad}M_{b,c}-\delta _{ac}M_{b,d}-\delta _{bd}M_{a,c}.
$$
The generators of $g=so(M)$ can be realized as follows. The Cartan
subalgebra ${\cal H}(g)$ is generated by
\begin{equation}
\label{sonorm-0}H_{j}=\left( -\frac i2\right) M_{2j-1,2j} ,\qquad j=1,...,N.
\end{equation}
For Cartan generators we shall also use the notation:
\begin{equation}
\label{sonorm-1}H_{j \pm \left( j+1\right) }=\left( -\frac i2\right) \left(
M_{2j-1,2j} \pm M_{2j+1,2j+2}\right).
\end{equation}
To the roots of $\Lambda(so(M))$ we attribute the generators
\begin{equation}
\label{sonorm-4}\left. 
\begin{array}{l}
E_{i+j}=\frac 12\left(
-M_{2i,2j}+iM_{2i,2j-1}+iM_{2i-1,2j}+M_{2i-1,2j-1}\right) ; \\ 
E_{i-j}=\frac 12\left(
-M_{2i,2j}-iM_{2i,2j-1}+iM_{2i-1,2j}-M_{2i-1,2j-1}\right) ; \\ 
E_{-i+j}=\frac 12\left(
+M_{2i,2j}-iM_{2i,2j-1}+iM_{2i-1,2j}+M_{2i-1,2j-1}\right) ; \\ 
E_{-i-j}=\frac 12\left(
+M_{2i,2j}+iM_{2i,2j-1}+iM_{2i-1,2j}-M_{2i-1,2j-1}\right) ; 
\end{array}
\right\} \qquad i<j; 
\end{equation}
and 
\begin{equation}
\label{sonorm-3} E_{ \pm k}=\frac 1{\sqrt{2}}\left( \pm
M_{2k,2N+1}-iM_{2k-1,2N+1}\right) , \qquad k\leq N, \quad {\rm for\ }%
so\left( 2N+1\right) . 
\end{equation}

The Borel subalgebras $B(g)$ are generated by the sets $\left\{ H_i, E_i,
E_{i \pm j} \right\} $ for $g=so(2N+1)$ and $\left\{ H_i, E_{i \pm j}
\right\} $ for $g=so(2N)$. To describe the chains of Frobenius subalgebras
we shall also need the alternative realization of these generators through
the ordinary matrix units. To get it let us renumerate the generators:
\begin{equation}
\left. 
\begin{array}{lcl}
A_{i,j} & \equiv & - E_{i-j}, \\ 
A_{i,2N+2-j} & \equiv & -E_{i+j}, \\ 
A_{i,N+1} & \equiv & -E_{i}, 
\end{array}
\right\} \qquad {\rm for} \quad so(2N+1) 
\end{equation}
and 
\begin{equation}
\left. 
\begin{array}{lcl}
A_{i,j} & \equiv & - E_{i-j}, \\ 
A_{i,2N+1-j} & \equiv & -E_{i+j},
\end{array}
\right\} \qquad {\rm for} \quad so(2N). 
\end{equation}

In these terms the Borel subalgebra $B(so(M))$ is spanned by the set $%
\left\{ A_{i,j} \mid i\leq j\right\} $ and  
we can also use the following matrix realization: 
\begin{equation}
\label{a-def}
\begin{array}{l}
H_i=\frac 12\left( 
{\cal E}_{i,i}-{\cal E}_{M+1-i,M+1-i}\right) , \\ A_{i,j}={\cal E}_{i,j}-%
{\cal E}_{M+1-j,M+1-i}.
\end{array}
\end{equation}

The canonical chains of twists (\ref{chain-ini}) for orthogonal simple Lie
algebras are based on the sequence of injections 
\begin{equation}
\label{soetower}
\begin{array}{l}
U(so(M))\supset U(so(M-4))\supset \ldots \supset U(so(M-4k))\supset \ldots
\\
\end{array}
\,
\end{equation}
Each link of such chains (see (\ref{chain-ini})) contains the jordanian
twist $\Phi _{{\cal J}_k}=$$\exp \{H_{\lambda _0^k} \otimes \sigma _{\lambda
_0^k}\}$ based on one of the long roots:
\begin{equation}
\label{ini-root}\lambda _0^k=e_1^k+e_2^k\quad ({\rm for \,\, both}\,\,D_N\
{\rm and}\ B_N)).
\end{equation}

The hyperplane $V_{\lambda _0^k}^{\perp }$ orthogonal to $\lambda _0^k$
coincides with the root space of the subalgebra $g_{\lambda _0^k}^{\perp
}=so\left( M-4\left( k+1\right) \right) \oplus so^{\left( k+1\right) }(3)$ .
For each subalgebra $g_{k+1}=so\left( M-4\left( k+1\right) \right) $ we can
again consider independently its root system $\Lambda^{k+1}$ and choose the
next initial root to be
$$
\lambda _0^{k+1}=e_1^{k+1}+e_2^{k+1}. 
$$
On each step we can fix two vector subrepresentations $d^{v\left( a\right) }$
in the restriction ( to${\rm \ }g_{k+1}$ ) of the adjoint representation $%
{\rm ad}_{g_k}$: 
$$
d_{g_{k+1}}^{v\left( a\right) }\subset {\rm ad}_{g_k\downarrow
g_{k+1}},\qquad a=1,2; 
$$
It is easy to check that the constituent roots form the weight diagrams for
these representations. The representation space for $d_{g_{k+1}}^{v\left(
a\right) }$ is spanned by the generators 
$$
\left\{ E_a,E_{a\pm l}\right\}\, {\rm with \, the \, roots}\, \left\{
e_a^k, \, e_a^k\pm e_l^k\in \pi _k\right\} 
$$
for $M-4k=2N+1$ and 
$$
\left\{E_{a\pm k}\right\}\,{\rm with \, the \, roots} \, \left\{ e_a^k\pm
e_l^k\in \pi _k\right\} 
$$
for $M-4k=2N$. In both cases $l=3,\ldots ,N$.

For the representations $d_{g_{k+1}}^{v\left( a\right) }$ and $%
d_{g_{k+1}}^{v\left( a\right) }\otimes d_{g_{k+1}}^{v\left( b\right) }$
the following scalar and tensor invariants, $I_{M-4k}^a$ and $%
I_{M-4k}^{a\otimes b}$ (with $a,b=1,2$ ), will be used in the construction
of twists and twisted coproducts:
\begin{equation}
\label{inv-1}
\begin{array}{l}
I_{2N+1}^a=\frac 12E_a^2+\sum_{l=3}^N\left( E_{a+l}E_{a-l}\right) , \\
I_{2N+1}^{a\otimes b}=E_a\otimes E_b+\sum_{l=3}^N\left( E_{a+l}\otimes
E_{b-l}+E_{a-l}\otimes E_{b+l}\right) , \\ 
I_{2N+1}^{a\land b}=E_a\land E_b+\sum_{l=3}^N\left( E_{a+l}\land
E_{b-l}+E_{a-l}\land E_{b+l}\right) , 
\end{array}
\end{equation}
\begin{equation}
\label{inv-2} 
\begin{array}{l}
I_{2N}^a=\sum_{l=3}^N\left( E_{a+l}E_{a-l}\right) , \\ 
I_{2N}^{a\otimes b}=\sum_{l=3}^N\left( E_{a+l}\otimes E_{b-l}+E_{a-l}\otimes
E_{b+l}\right) , \\ 
I_{2N}^{a\land b}=\sum_{l=3}^N\left( E_{a+l}\land E_{b-l}+E_{a-l}\land
E_{b+l}\right) . 
\end{array}
\end{equation}

The set of initial roots defines a natural gradation in the root space of
the subalgebra $N^{+}\left(so \left( M \right)\right) $ $\subset $ $so\left(
M\right) $ : 
\begin{equation}
\label{gradation}\Lambda \left( N^+(so\left( M\right) )\right)
=\bigcup_{k=0}^{p_{\max} }\left( \lambda _0^k\cup \pi _k\right) , 
\end{equation}
where $p_{\max }=\left[ M/4\right] +\left[ \left( M+1\right) /4\right] $.

The inverse of the map defined by the classical $r$-matrix is the Frobenius
bilinear form. Let us study the Frobenius subalgebras in $B(so(M))$. 
\newtheorem{proposition}{Proposition} 
\begin{proposition}
Let ${\bf L}$ be a semi-direct sum of a subalgebra  ${\bf S}$ and a
commutative ideal ${\bf N}$.
 Then ${\bf L}$ is Frobenius if and only if the following
conditions hold:
\begin{description}
\item[i)]
${\bf L}$ acts transitively on the space ${\bf N} ^*$ with the generic
point $A ^*$;
\item[ii)]
the stationary subalgebra ${\bf S}_{A^*}
= \left\{ s \in  {\bf S}: A ^*([s,x])
 = 0, {\rm for \,  any } \, x \in {\bf N}  \right\}$ is Frobenius with
a Frobenius homomorphism $f_0 :{\bf S}_{A ^*} \longrightarrow {\bf C} $.
\end{description}

Moreover, in this case $f = f_0 \oplus A ^*$ is a Frobenius homomorphism
for ${\bf L}$.
\end{proposition}
This statement can be obtained as a consequence of the Proposition 1 and the
Remark after it in \cite{STO2}. Here is how it can be used in the case of
orthogonal simple Lie algebras. \newtheorem{lemma}{Lemma} 
\begin{lemma}
Let ${\bf L}_1 \subset B(so(M)) $ be a subalgebra generated by the set
\\ $ \left\{ H_1, H_2, A_{i,j}, i = 1,2
\right\}$. Then  ${\bf L}_1$ is Frobenius.
\end{lemma}
{\bf Proof}. In ${\bf L}_1$ the following subalgebras can be fixed ${\bf N}%
_1=\left\{ A_{1,j}\right\} ,{\bf S}_1=\left\{ H_1,H_2,A_{2,j}\right\} $. The
generators of the dual space ${\bf N}_1^{*}$ can be identified with the
matrices $\left\{ A_{i,1}\right\} $ defined according to the rule (\ref
{a-def}) and connected with $A_{1,j}$ through the bilinear form $<A,B>={\rm %
tr}(A,B)$ or by a Killing form in the general setting. Since dim${\bf N}_1=$
dim${\bf S}_1$ it suffices to find a point $A^{*}\in {\bf N}_1^{*}$ such
that ${\bf S}_{A^{*}}=0$. One can check directly that $A_0^{*}=%
\sum_{i=2}^{n-1}A_{i1}$ satisfies this condition. $\spadesuit $

This point $A_0^*$ is not unique. If $G({\bf S}_1)$ is the subgroup of $%
SO(M) $ corresponding to the algebra ${\bf S}_1 $ then for any $g \in G({\bf %
S}_1) $ the point ${\rm Ad}^*(g)(A_0^* ) = (A_0^*)^g $ satisfies the
condition ${\bf S}_{ (A_0^*)^g}=0 $ since ${\bf S}_{ (A_0^*)^g}= g^{-1}{\bf S%
}_{A_0^*} g = 0$. For our purposes it is convenient to choose $A_0^*
=A_{2,1} + A_{M-1,1}$ (One can check that this point satisfies the
conditions of the Proposition 1.)

\begin{lemma}
Let ${\bf L}_{K,M}$ be a subalgebra of $B(so(M))$ generated by the set
$$
\left\{ H_i, A_{i,j} \mid i= 1, \dots ,2K; j=1, \dots M; i < j ; 2K \leq [M/2]
\right\}.
$$
Then ${\bf L}_{K,M}$ is Frobenius.
\end{lemma}

\noindent {\bf Proof}. The algebra ${\bf L}_{K,M}$ has the structure of a
semidirect sum: 
$$
{\bf L}_{K,M} = {\bf S}_{K,M} \dot{\oplus} {\bf N}_{K,M},
$$
where 
$$
{\bf N}_{K,M} = \left\{ \left\{ H_i \mid i= 1, \dots ,2K; \right\} 
\left\{A_{i,j} \mid i= 2, \dots ,2K; j=1, 
\dots M; i < j \right\} \right\}. 
$$
Evidently ${\bf L}_{K,M}$ acts transitively on ${\bf N}_{K,M}^*$ with the
generic point $A_1^* = A_{2,1} + A_{M-1,1} $. One can easily check that 
$$
\left( {\bf S}_{K,M} \right) _{ (A_0^*)} = \left\{ H_i, A_{i,j} 
\mid i= 3, \dots ,2K; i < j \right\}. 
$$
Thus 
$$
\left( {\bf S}_{K,M} \right) _{ (A_0^*)} \cong {\bf L}_{K-1,M-4} \subset
B(so(M-4)). 
$$
Obvious induction shows that ${\bf L}_{K,M}$ is Frobenius due to the
Proposition 1. $\spadesuit$

The algebra ${\bf L}_{K,M}$ has a nontrivial second cohomology group 
$H^2({\bf L}_{K,M})$. 
The latter contains $\Lambda ^2\left( H_{K,M}^{*}\right) $
where $H_{K,M}^{*}$ is the space dual to the Cartan subalgebra in 
${\bf L}_{K,M}$, $H_{K,M}^{*}= 
{\bf L}_{K,M}\bigcap {\cal H}(so(M))\subset B(so(M))$.
It is easy to see that all bilinear forms $H_i^{*}\bigwedge H_j^{*}$ are
2-cocycles and not coboundaries. Here $H_i^{*}\in H_{K,M}^{*}$. Consequently 
$A_{K,M}^{*}+\zeta _{ij}H_i^{*}\bigwedge H_j^{*}$ are the nondegenerate
2-cocycles on ${\bf L}_{K,M}$. The map $A_{K,M}^{*}$ is a Frobenius
homomorphism, $A_{K,M}^{*}:{\bf L}_{K,M}\longrightarrow {\bf C}$, because $%
A_{K,M}^{*}\left( [x,y]\right) $ is a nondegenerate bilinear form on ${\bf L}%
_{K,M}$. The induction procedure shows that $A_{K,M}^{*}$ can be chosen in
the following form: 
\begin{equation}
\label{gen-form}A_{K,M}^{*}=\left( A_{21}+A_{43}+A_{65}\dots \right) +\left(
A_{M-1,1}+A_{M-3,3}+A_{M-5,5}\dots \right) .
\end{equation}

In the case of the orthogonal simple Lie algebras chains of twists (\ref
{chain-ini}) introduced in \cite{KLO} refer to the Frobenius subalgebras
(contained in the corresponding $B(so\left( M\right) )$) with the coboundary
nondegenerate bilinear forms \cite{STO},
\begin{equation}
\label{ini-fro}\omega _p^{+}=\sum_{k=0}^p\gamma _k\left( E_{1+2}^{\left(
k\right) }\right) ^{*}\left( \left[ \, , \, \right] \right) ,
\end{equation}
where $p$ is the number of links in the chain ${\cal F}_{{\cal B}_{0\prec
p}} $ and the parameters $\gamma _k=0,1$ indicate that we describe here the
set of forms. It is obvious that the forms (\ref{ini-fro}) ignore the
subspaces of $so^{(k)}(3) $ subalgebras that appear
on each step of the sequence
of injections in (\ref{gradation}). According to the formula (\ref{gen-form}%
) to describe the maximal Frobenius subalgebras in $B(so(M))$ the following
form must be considered:
\begin{equation}
\label{new-fro}\omega _p^{\pm }=\sum_{k=0}^p\left( \gamma _k\left(
E_{1+2}^k\right) ^{*}-\delta _k\left( E_{1-2}^k\right) ^{*}\right) \left(
\left[ \, , \, \right] \right) .
\end{equation}
Here both parameters are discrete: $\gamma _k,\delta _k=0,1$ . (Notice that
this does not lead to the undesirable terms in the corresponding carrier
space because in the Borel subalgebra $B(so\left( M\right) )$ (fixed by the
choice of $\left\{ \lambda _0^k\right\} $) there are no constituent roots
for $e_1^k-e_2^k$. The form $\omega _p^{\pm }$ is considered on $B(so\left(
M\right) )$. In terms of the integral root system $\Lambda(so(M))$ (not
split in the subsystems $\Lambda^{(k)}$) the generators $E_{1+2}^{\left(
k\right) }$ and $E_{1-2}^{\left( k\right) }$ form the sequence : 
$$
\begin{array}{c}
\left\{ E_{1+2},E_{1-2},E_{3+4},E_{3-4},\dots ,E_{(2 p+1) +(2 p+2) }, 
E_{(2p+1) -(2p+2) }\right\} \approx \\
\approx \left\{ A_{1,M-1},A_{1,2},A_{3,M-3},A_{3,4},\dots A_{2p+1,M-(2p+1)
},A_{2p+1,2p+2 },\right\}. 
\end{array}
$$

Thus we come to the conclusion that there must be two sets of chains of
twists for the orthogonal algebras corresponding to the two sets of the
coboundary forms (\ref{ini-fro}) and (\ref{new-fro}). The first set is the
canonical chain of twists (\ref{chain-ini}), whose twisting element can be
rewritten in terms of invariants (\ref{inv-1}),(\ref{inv-2}):
\begin{equation}
\label{chain} 
\begin{array}{l}
{\cal F}_{{\cal B}_{0\prec p}}=\exp \left\{ I_{M-4p}^{1\otimes 2}\left(
1\otimes e^{-\frac 12\sigma _{1+2}^p}\right) \right\} \cdot \exp
\{H_{1+2}^{\left( p\right) }\otimes \sigma _{1+2}^p\}\,\cdot \\ \exp \left\{
I_{M-4\left( p-1\right) }^{1\otimes 2}\left( 1\otimes e^{-\frac 12\sigma
_{1+2}^{p-1}}\right) \right\} \cdot \exp \{H_{1+2}^{\left( p-1\right)
}\otimes \sigma _{1+2}^{p-1}\}\,\cdot \\ 
\ldots \\
\exp \left\{ I_M^{1\otimes 2}\left( 1\otimes e^{-\frac 12\sigma
_{1+2}^0}\right) \right\} \cdot \exp \{H_{1+2}^{\left( 0\right) }\otimes
\sigma _{1+2}^0\}\, \\ 
=\prod_{k=p}^0\exp \left\{ I_{M-4k}^{1\otimes 2}\left( 1\otimes e^{-\frac
12\sigma _{1+2}^k}\right) \right\} \cdot \exp \{H_{1+2}^{\left( k\right)
}\otimes \sigma _{1+2}^k\} 
\end{array}
\end{equation}
(here $\sigma _{\lambda _0^k}$ are denoted by $\sigma _{1+2}^k=\ln \left(
1+E_{1+2}^{\left( k\right) }\right) $ according to (\ref{ini-root}) and for
simplicity we put $\gamma _l=1$ ).

When dealing with the forms $\omega _p^{\pm }$ the problem is that in the
process of twisting by a chain (\ref{chain}) the costructure of the
subalgebras $so^{\left( k\right) }(3)$ is considerably changed and the twist
equations (\ref{twist-eq}) become extremely difficult to solve.

\section{Construction of the full chains of twists}

According to the general structure of a chain of twists \cite{KLO} we can
study its links separately. Let us assume that we have constructed the $k-1$
links of a chain and found the matreshka effect. This means that after the
chain twisting with $k-1$ links we get the subalgebra $g^{\left( k\right)
}=so\left(M-4k\right) $ with primitive generators. We shall show that it is
possible to construct the next link of the chain so that the twist
will correspond to the enlarged form $\omega _p^{\pm }$ (see (\ref
{new-fro})). To start the construction of the $k$-th link we have to choose
the initial root $\lambda _0^k$ (as in (\ref{ini-root})) and the subalgebra $%
{\bf L}_{K,M-4k}$ described in Lemma 2 (with $2K=[M/2]$). First we apply the
following jordanian twist to the subalgebra ${\bf L}_{K,M-4k}$:
\begin{equation}
\label{jord-ini}\Phi _{{\cal J}_k}=\exp \left( H_{1+2}^k\otimes \sigma
_{1+2}^k\right). 
\end{equation}
This results in the following deformed coproducts: 
$$
\begin{array}{l}
\Delta _{
{\cal J}_k}\left( H_{1+2}^k\right) =H_{1+2}^k\otimes e^{-\sigma
_{1+2}^k}+1\otimes H_{1+2}^k, \\ \Delta _{
{\cal J}_k}\left( E_{1+2}^k\right) =E_{1+2}^k\otimes e^{\sigma
_{1+2}^k}+1\otimes E_{1+2}^k, \\ \Delta _{
{\cal J}_k}\left( E_{a\pm l}^k\right) =E_{a\pm l}^k\otimes e^{\frac 12\sigma
_{1+2}^k}+1\otimes E_{a\pm l}^k, \\ \quad l=3,\ldots ,N;\quad a=1,2; 
\end{array}
$$
For $M=2N+1$ we also get 
$$
\Delta _{{\cal J}_k}\left( E_a^k\right) =E_a^k\otimes e^{\frac 12\sigma
_{1+2}^k}+1\otimes E_a^k,
$$
Notice that the generators $\left\{ H_{1-2}^k,E_{1-2}^k,E_{2-1}^k\right\} $
remain primitive.

The second twisting factor must be the full canonical extension \cite{KLM}
for the jordanian twist $\Phi _{{\cal J}_k}$ (\ref{jord-ini}): 
\begin{equation}
\label{ext-ini}\Phi _{{\cal _E}_k}=\exp \left( I_{M-4k}^{1\otimes 2}\left(
1\otimes e^{-\frac 12\sigma _{1+2}^k}\right) \right) , 
\end{equation}
The successive application of these two factors performs the extended
jordanian twisting by the element $\Phi _{{\cal _E}_k}\Phi _{{\cal J}_k}$
that leads to the following costructure in ${\bf L}_{K,M-4k}$: 
$$
\begin{array}{l}
\Delta _{
{\cal E}_k{\cal J}_k}\left( H_{1+2}^k\right) =H_{1+2}^k\otimes e^{-\sigma
_{1+2}^k}+1\otimes H_{1+2}^k-\left( 1\otimes e^{-\frac 32\sigma
_{1+2}^k}\right) I_{M-4k}^{1\otimes 2}, \\ \Delta _{
{\cal E}_k{\cal J}_k}\left( H_{1-2}^k\right) =H_{1-2}^k\otimes 1+1\otimes
H_{1-2}^k, \\ \Delta _{
{\cal E}_k{\cal J}_k}\left( E_{1+2}^k\right) =E_{1+2}^k\otimes e^{\sigma
_{1+2}^k}+1\otimes E_{1+2}^k, \\ 
\Delta _{{\cal E}_k{\cal J}_k}\left( E_{1\pm k}^k\right) 
=E_{1 \pm k}^k\otimes e^{-\frac 12\sigma _{1+2}^k}+
1\otimes E_{1 \pm k}^k, \\ 
\Delta _{{\cal E}_k{\cal J}_k}\left( E_{2\pm k}^k\right) 
=E_{2 \pm k}^k\otimes e^{\frac 12\sigma
_{1+2}^k}+e^{\sigma _{1+2}^k}\otimes E_{2 \pm k}^k, 
\end{array}
$$

$$
\begin{array}{l}
\Delta _{
{\cal E}_k{\cal J}_k}\left( E_{1-2}^k\right) = \\ E_{1-2}^k\otimes
1+1\otimes E_{1-2}^k+\left( 1\otimes e^{-\frac 12\sigma _{1+2}^k}\right)
I_{M-4k}^{1\otimes 1}+I_{M-4k}^1\otimes \left( e^{-\sigma _{1+2}^k}-1\right)
, \\ 
\Delta _{
{\cal E}_k{\cal J}_k}\left( E_{2-1}^k\right) = \\ E_{2-1}^k\otimes
1+1\otimes E_{2-1}^k+\left( e^{\sigma _{1+2}^k}-1\right) \otimes
I_{M-4k}^2e^{-\sigma _{1+2}^k}+\left( 1\otimes e^{-\frac 12\sigma
_{1+2}^k}\right) I_{M-4k}^{2\otimes 2}. 
\end{array}
$$
And in the case of $M=2N+1$ for the short root generators we get: 
$$
\begin{array}{l}
\Delta _{
{\cal E}_k{\cal J}_k}\left( E_1^k\right) =E_1^k\otimes e^{-\frac 12\sigma
_{1+2}^k}+1\otimes E_1^k, \\ \Delta _{{\cal E}_k{\cal J}_k}\left(
E_2^k\right) =E_2^k\otimes e^{\frac 12\sigma _{1+2}^k}+e^{\sigma
_{1+2}^k}\otimes E_2^k. 
\end{array}
$$

It would be necessary to have the coproducts for some of the invariants (see
(\ref{inv-1}) and (\ref{inv-2})) 
\begin{equation}
\label{co-i1} 
\begin{array}{l}
\Delta _{{\cal E}_k{\cal J}_k}\left( I_{M-4k}^1\right) = \\ 
I_{M-4k}^1\otimes e^{-\sigma _{1+2}^k}+1\otimes I_{M-4k}^1+ 
I_{M-4k}^{1\otimes 1}\left(1\otimes e^{-\frac 12\sigma _{1+2}^k}\right) , 
\end{array}
\end{equation}
\begin{equation}
\label{co-i2} 
\begin{array}{l}
\Delta _{
{\cal E}_k{\cal J}_k}\left( I_{M-4k}^2e^{-\sigma _{1+2}^k}\right) = \\ 
I_{M-4k}^2e^{-\sigma _{1+2}^k}\otimes 1+e^{\sigma _{1+2}^k}\otimes
I_{M-4k}^2e^{-\sigma _{1+2}^k}+I_{M-4k}^{2\otimes 2}\left( 1\otimes
e^{-\frac 12\sigma _{1+2}^k}\right) . 
\end{array}
\end{equation}

We have two generators of ${\bf L}_{K,M-4k}$ that are not yet incorporated
in the carrier subalgebra of the twist: $H_{1-2}^k$ and $E_{1-2}^k$. The
coproduct of the latter is deformed. So the canonical jordanian factor
cannot be used here. In \cite{KL} it was indicated that the reason of the
nonprimitivity of the coproduct $\Delta_{{\cal E}_k{\cal J}_k}\left(
E_{1-2}^k \right)$ is that the generator $E_{1-2}^k$ belongs to the long
series of the initial root $\lambda _0^k = e_1^k + e_2^k$. It was shown
there that in such a case the deformed carrier subspace must exist with
primitive basic elements. In our situation such ''deformed'' generators must
have the form 
\begin{equation}
\label{g1-2} 
\begin{array}{l}
G_{1-2}^{k}=E_{1-2}^k-I_{M-4k}^1, \\ 
G_{2-1}^{k}=E_{2-1}^k-I_{M-4k}^2e^{-\sigma _{1+2}^k}. 
\end{array}
\end{equation}
Using the coproducts (\ref{co-i1}, \ref{co-i2}) it is easy to check that
both $G_{1-2}^{k}$ and $G_{2-1}^{k}$ are primitive,
$$
\begin{array}{l}
\Delta _{
{\cal E}_k{\cal J}_k}\left( G_{1-2}^{k}\right) =\Delta _{{\cal E}_k{\cal J}%
_k}\left( E_{1-2}^k\right) -\Delta _{{\cal E}_k{\cal J}_k}\left(
I_{M-4k}^1\right) \\ =G_{1-2}^{k}\otimes 1+1\otimes G_{1-2}^{k}, \\
\Delta _{{\cal E}_k{\cal J}_k}\left( G_{2-1}^{k}\right) 
=\Delta _{{\cal E}_k{\cal J}_k}\left( E_{2-1}^k\right) 
-\Delta _{{\cal E}_{k}{\cal J}_{k}}\left(
I_{M-4k}^2\right) \left( e^{-\sigma _{1+2}^k}\otimes e^{-\sigma
_{1+2}^k}\right) \\ =G_{2-1}^{k}\otimes 1+1\otimes G_{2-1}^{k}.
\end{array}
$$

Together with $H_{1-2}^k$ the elements (\ref{g1-2}) generate a 3-dimensional
space $V_G^{k}$ of primitive elements in the algebra $U_{{\cal E}_k{\cal J}%
_k}\left( so(M-4k)\right) $. Both $G_{1-2}^{k}$ and $G_{2-1}^{k}$ commute
with $U\left( so\left( M-4(k+1)\right) \right) $ as well as $H_{1-2}^k$
whose dual vector is orthogonal to the roots of $so\left( M-4(k+1) \right) $.

The subspace $V_G^{k}$ spanned by $\left\{
H_{1-2}^k,G_{1-2}^{k},G_{2-1}^{k}\right\} $ is algebraically closed:
$$
\begin{array}{l}
\left[ H_{1-2}^k,G_{1-2}^{k}\right] =G_{1-2}^{k}, \\ 
\left[ H_{1-2}^k,G_{2-1}^{k}\right] =-G_{2-1}^{k}, \\
\left[ G_{1-2}^{k},G_{2-1}^{k}\right] =2H_{1-2}^k.
\end{array}
$$
Let us denote this algebra by $so_G^{\left( k \right) }(3)$. Clearly it is
primitive, commutes with $U\left( so\left( M-4\left( k+1\right) \right)
\right)$ and is realized on a deformed subspace. (This space is not
orthogonal to $H_{1+2}^k$. Moreover, $G_{2-1}^{k},G_{1-2}^{k}$ are not any
longer eigenvectors of ${\rm ad}_{H_{1+2}^k}$.)

Another subalgebra which remains primitive after the composition $\Phi _{%
{\cal E}_k}\Phi _{{\cal J}_k}$ of twists (\ref{jord-ini}) and (\ref{ext-ini}%
) is $so\left( M-4\left( k+1\right) \right) $ (due to the matreshka effect).
We come to conclusion that the twisted $U_{{\cal E}_k{\cal J}_k}\left(
so\left( M-4k\right) \right) $ contains the primitive subalgebra $g_{\lambda
_0^k}^{\perp }$ 
$$
U_{{\cal E}_k{\cal J}_k}\left( so\left( M-4k\right) \right) \supset
g_{\lambda _0^k}^{\perp }=so\left( M-4\left( k+1\right) \right) \oplus
so_G^{\left( k\right) }(3). 
$$
Its Borel subalgebra is ${\bf L}_{K,M-4(k+1)}\oplus B(so_G^{\left( k\right)}
(3))$ and it is Frobenius (see Section 2).

Remember that the subalgebra $g_{\lambda _0^k}^{\perp }$ has a structure of
direct sum. Further, twisting by the next factors (such as $\Phi _{{\cal E}%
_{k+s}} \Phi _{{\cal J}_{k+s}}$) can not affect the primitive subalgebra $%
so_G^{\left( k\right) }(3)$. Each step produces (in the corresponding $%
g_{\lambda _0^k}^{\perp }$) the additional subalgebra $so_G^{\left( k\right)
}(3), \quad k=1,\ldots ,p$. The primitive subalgebras that can be found in
an orthogonal algebra after the chain twisting (\ref{chain}) with $p$ links
contain not only $so\left( M-4\left( p+1\right) \right) $ but also a direct
sum of $p$ copies of $so_G(3)$: 
$$
U_{{\cal B}_{0\prec p}}\left( so\left( M\right) \right) \supset {\cal D}=
\oplus_{k=1}^p \, so_G^{\left( k \right) }(3). 
$$

The main consequence is that in the twisted $U_{{\cal B}_{0\prec p}}\left(
so\left( M\right) \right) $ one can perform further twist deformations with
the carrier subalgebra in ${\cal D}$ . The most interesting among them are
the jordanian twists defined by 
\begin{equation}
\label{g-twist}\Phi _{{\cal J}_k}^G=\exp \left( H_{1-2}^k\otimes \sigma
_{_G}^k\right) 
\end{equation}
that can be attributed to any number of copies $so_G(3)$ . Here $\sigma
_{_G}^k = \ln \left( 1 \right. + \left. G_{1-2}^{k}\right) .$ Thus in the
general expression for the twisting element (\ref{chain}) one can insert in
the appropriate $k$ places the additional factors which are the jordanian
twisting elements on the deformed carrier spaces. This means that we can
perform a substitution 
$$
\begin{array}{l}
\Phi _{
{\cal E}_k}\Phi _{{\cal J}_k}\Rightarrow \Phi _{{\cal J}_k}^G\Phi _{{\cal E}%
_k}\Phi _{{\cal J}_k}=\Phi _{{\cal G}_k} \\ \exp \left\{ I_{M-4k}^{1\otimes
2}\left( 1\otimes e^{-\frac 12\sigma _{1+2}^k}\right) \right\} \cdot \exp
\{H_{1+2}^k\otimes \sigma _{1+2}^k\}\Rightarrow \\
\exp \left( H_{1-2}^k\otimes \sigma _{_G}^k\right) \cdot \exp \left\{
I_{M-4k}^{1\otimes 2}\left( 1\otimes e^{-\frac 12\sigma _{1+2}^k}\right)
\right\} \cdot \exp \left( H_{1+2}^k\otimes \sigma _{1+2}^k\right) 
\end{array}
$$
Thus {\em the full chain} has the following form 
\begin{equation}
\label{chain-new} 
\begin{array}{l}
{\cal F}_{{\cal G}_{0\prec p}}=\prod_{k=p}^0\Phi _{{\cal G}_k}= \\ 
\prod_{k=p}^0\exp \left( H_{1-2}^k\otimes \sigma _{_G}^k\right) \cdot \exp
\left\{ I_{M-4k}^{1\otimes 2}\left( 1\otimes e^{-\frac 12\sigma
_{1+2}^k}\right) \right\} \cdot \exp \{H_{1+2}^k\otimes \sigma _{1+2}^k\}. 
\end{array}
\end{equation}

This result means that we have constructed the explicit quantizations
with a triangular universal $R$-matrix
$$
{\cal R}_{{\cal G}_{0\prec p}}=\left( {\cal F}_{{\cal G}_{0\prec p}}\right)
_{21}\left( {\cal F}_{{\cal G}_{0\prec p}}\right) ^{-1} 
$$
for the following set of classical $r$-matrices:%
$$
r_{{\cal G}_{0\prec p}}=\sum_{k=0}^p\eta _k\left( H_{1+2}^k\land
E_{1+2}^k+\xi _kH_{1-2}^k\land E_{1-2}^k+I_{M-4k}^{1\land 2}\right).  
$$
Here all the parameters are independent and continuous. Elementary
computations show that these full chains (\ref{chain-new}) correspond to the
coboundary forms (\ref{new-fro}). To illustrate these quantizations we
present in Appendix the matrix ${\cal R}_{{\cal G}_{0\prec p}}$ for the
algebra $so(5)$ in the defining representation.

In Section 2 we proved that adding 
$\zeta_{ij} H_{i}^* \bigwedge H_{j}^*$ to
the forms of type (\ref{new-fro}) we obtain new non-degenerate 2-cocycles,
which are not coboundaries. We can also construct the corresponding twists
for these modified cocycles: 
\begin{equation}
\label{new-fro-nc} \omega _p^{\pm }=\sum_{k=0}^p\left( \gamma _k\left(
E_{1+2}^k\right) ^{*}-\delta _k\left( E_{1-2}^k\right) ^{*}\right) \left(
\left[ \, , \, \right] \right) + \sum_{i,j=0; \, i \neq j}^p \zeta_{ij}
H_{i}^* \bigwedge H_{j}^* .
\end{equation}
Notice that the subalgebras $so_G^{\left( k\right) }(3)$ commute not only
with $so\left( M-4(k + 1)\right)$ but also with any 
$\left\{ E_{1 + 2}^{(s)} \,| \,s \leq k \right\}$. 
This means that after having twisted $U \left(
so\left( M\right) \right)$ by the chain (\ref{chain-new}) we obtain $p+1$
pairs of commuting primitive elements 
$\left\{ \sigma^k_{1+2}, \sigma^k_{G}\,|\, k = 0, \dots , p \right\} $. 
Therefore we can apply the Reshetikhin twist 
\begin{equation}
\label{r-twist}\Phi _{{\cal R}}=\exp \left( \zeta_{ij} \sigma_i \otimes
\sigma_j \right), \quad \sigma_i \in \left\{ \sigma^k_{1+2}, 
\sigma^k_{G}\, |\, k = 0, \dots , p \right\}. 
\end{equation}
to the algebra $U_{{\cal G}_{0\prec p}}\left( so\left( M\right) \right)$.

Thus the element 
$$
\Phi _{{\cal R}}{\cal F}_{{\cal G}_{0\prec p}} 
$$
defines also a twist for $U \left( so\left( M\right) \right)$. It leads to
the deformed  Hopf algebra $U_{{\cal R G}_{0\prec p}}$$\left( so\left(
M\right)\right)$ with the universal element
$$
{\cal R}_{{\cal R G}_{0\prec p}}= \left(\Phi _{{\cal R}} {\cal F}_{{\cal G}%
_{0\prec p}}\right) _{21}\left( {\cal F}_{{\cal G}_{0\prec p}}\right)
^{-1}\left( \Phi _{{\cal R}}\right) ^{-1}
$$
and the classical $r$-matrix
$$
\begin{array}{lcl}
r_{{\cal R G}_{0\prec p}} & = & \sum_{k=0}^p\eta _k\left( H_{1+2}^k\land
E_{1+2}^k+\xi _kH_{1-2}^k\land E_{1-2}^k+I_{M-4k}^{1\land 2}\right) \\
&  & + \sum_{i,j=0; \, i \neq j}^p \zeta_{ij} E_{s}^i \land E_{t}^j; \\
&  & \rule{0cm}{1cm} 
E^k_{s},E^k_{t} \in \left\{ E^k_{1+2}, E^k_{1-2}\, |\, k =
0, \dots , p \right\} .
\end{array}
$$

The dimensions of the nilpotent subalgebras $N^{+}\left( so\left( M\right)
\right) $ in the sequence $g_{\lambda _0^p}^{\perp }\subset g_{\lambda
_0^{p-1}}^{\perp }\subset \ldots \subset g_{\lambda _0^0}^{\perp }\subset g$
are subject to the following simple relation:
\begin{equation}
\label{dim}\dim \left( N^{+}\left( so\left( M\right) \right) \right) -\dim
\left( N^{+}\left( so\left( M-4\right) \right) \right) =2\left( \dim
d_{so\left( M-4\right) }^v+1\right) .
\end{equation}
Taking this into account we conclude (from the formula (\ref{gradation}))
that the chains (\ref{chain-new}) are full. Furthermore, this means that for
$p=p_{\max }=\left[ M/4\right] +\left[ \left( M+1\right) /4\right] $ the
corresponding carrier spaces contain all the generators of the nilpotent
subalgebra $N^{+}\left( so\left( M\right) \right) $. When $M$ is even-odd
one can always find in $so\left( M\right) $ one independent Cartan generator
which cannot be included in the carrier subalgebra of a chain. When $M$ is
even-even or odd the total number of jordanian twists in a maximal full
chain ${\cal F}_{{\cal G}_{0\prec p\max }}$ is equal to the rank of $%
so\left( M\right) $. Thus in the latter case the carrier subalgebra is equal
to the Borel subalgebra. 

\section{Conclusions}

The family of explicit twisting elements was constructed for the universal
enveloping algebras ${\cal A}=U(so(M))$ (series $B$, and $D$) with full
nilpotent subalgebras $N^{+}\left( so\left( M\right) \right) $ included in
the corresponding carrier spaces.

There is a variety of applications for explicitly known twisting elements $%
{\cal F}$. Using a particular (e.g. fundamental) representation for one of
the factors of ${\cal A}\otimes {\cal A}$ we get from the universal $R$%
-matrix the $L$-operator of the FRT-formalism and this results in explicit
relations among the generators of the original universal enveloping algebra
and the FRT-generators of the twisted one.

Twisting of the coalgebra in ${\cal A}$ induces changes in Clebsch-Gordan
coefficients of bases in the tensor products of irreducible representations $%
c_V\otimes d_W$. The evaluation of these coefficients is given by the direct
action of the matrix $F$ ( $F$ is the value of the twisting element in the
corresponding representation: $F=c_V\otimes d_W\left( {\cal F}\right) $ ) on
the original CG coefficients \cite{KS}.

Due to the embedding of the simple Lie algebras $g$ into the corresponding
Yangians (as Hopf subalgebras) $U(g)\subset {\cal Y}(g)$ \cite{DRQG} the
Yangian $R$-matrix $R_{{\cal Y}}$
can be twisted by the same ${\cal F}$ defined for $%
g$ \cite{KST}. As a result for the case of orthogonal algebra $g=so(M)$ the
$R$-matrix (in the defining representation $d\subset {\rm %
Mat}(M,{\bf C})\otimes {\rm Mat}(M,{\bf C})$ ) will be changed: 
$$
\begin{array}{c}
ud\left( 1\otimes 1\right) + 
{\cal P}-\frac u{u-1+M/2}{\cal K\quad \longrightarrow } \\ \quad ud\left( 
{\cal F}_{21}{\cal F}^{-1}\right) +{\cal P}-\frac u{u-1+M/2}d\left( {\cal F}%
_{21}\right) {\cal K} d \left( {\cal F}^{-1}\right) 
\end{array}
$$
(Here $u$ is a spectral parameter and the operator ${\cal K}$ is obtained
from the permutation
${\cal P}$ by transposing its first tensor factor.) For the canonical
chains ${\cal F}={\cal F}_{{\cal B}_{0\prec p}}$ the deformed solutions of
YBE were given in the explicit form in \cite{L}. Similarly to the case of
canonical chains (\ref{chain-ini}) 
the twists ${\cal F}_{{\cal G}_{0\prec p}}$ produce the
sets of deformed Yangians and the new integrable hamiltonians 
(cf. the $sl(2)$-case \cite{KS}).

\section{Acknowledgements}

This work was supported in part by the Russian Foundation for Basic Research
under the grants N 99-01-00101 (P.P.K.) and N 00-01-00500 (V.D.L.), the
Royal Swedish Academy of Sciences under the program ``Cooperation between
Sweden and former USSR'', the INTAS grant N 99-01459 (P.P.K.).

\section{Appendix}

Here we take $g = so(5)$. This is the simplest case where the full chain
differs nontrivially from the canonical one (\ref{chain-ini}) \cite{KLO}
and the deformed carrier space is
used to construct the twist:
\begin{equation}
\label{chain-so5} 
\begin{array}{l}
{\cal F}_{{\cal G}_{0\prec p}}= \\ \exp \left( H_{1-2}^0\otimes \sigma
_{_G}^0\right) \cdot \exp \left\{ I_{5}^{1\otimes 2}\left( 1\otimes
e^{-\frac 12\sigma _{1+2}^0}\right) \right\} \cdot \exp \{H_{1+2}^0\otimes
\sigma _{1+2}^0\}. 
\end{array}
\end{equation}
Consider now the corresponding $R$-matrix ${\cal R}_{{\cal G}_{0\prec
p}}=\left( {\cal F}_{{\cal G}_{0\prec p}}\right) _{21}\left( {\cal F}_{{\cal %
G}_{0\prec p}}\right) ^{-1} $ in the defining 5-dimensional representation $%
d $ of $so(5)$. This means that we use the following matrix realization for
the generators of $B(so(5))$:
\begin{equation}
\label{def} 
\begin{array}{l}
H_{i} = 
\frac{1}{2}\left( {\cal E}_{i,i} - {\cal E}_{M+1-i,M+1-i} \right), \\ 
A_{i,j} = 
{\cal E}_{i,j} - {\cal E}_{M+1-j,M+1-i}; \\ \rule{0cm}{1cm} i,j = 1, \dots ,
5. 
\end{array}
\end{equation}
As a result we get the solution of the matrix Yang-Baxter equation that can
be written in terms of tensor products of $5 \times 5$ matrix units 
${\cal E}_{i,j}$: 
$$
\label{rn-so5} 
\begin{array}{l}
d\left(
{\cal R}_{{\cal G}_{0\prec p}} \right)= R_{{\cal G}_0}= \\ (
{\cal E}_{4,5} - {\cal E}_{1,2}) \otimes \left( \frac{1}{2}( {\cal E}_{1,1}- 
{\cal E}_{2,2}+ {\cal E}_{4,4}- {\cal E}_{5,5}) + \frac{1}{4}( 
{\cal E}_{4,5}- {\cal E}_{1,2} ) 
- \frac{1}{8} {\cal E}_{1,5} \right) + \\ ({\cal E}_{4,5} + {\cal E}_{1,2}) 
\otimes \left( - \frac{1}{4}( {\cal E}_{1,4}+ {\cal E}_{2,5}) - 
\frac{1}{8} {\cal E}_{1,5} \right) + \\ ( 
{\cal E}_{3,5} - {\cal E}_{1,3}) \otimes 
\left( ( {\cal E}_{2,3}- {\cal E}_{3,4}) \right) + \\ ( 
{\cal E}_{3,5} + {\cal E}_{1,3}) \otimes \left( - \frac{1}{2} 
( {\cal E}_{3,5}+ {\cal E}_{1,3}) \right) + \\ 
({\cal E}_{3,4} - {\cal E}_{2,3} ) \otimes 
\left( ( {\cal E}_{3,5}- {\cal E}_{1,3}) \right) + \\ ( 
{\cal E}_{2,5} - {\cal E}_{1,4} ) \otimes \left( \frac{1}{4}( {\cal E}%
_{2,5}- {\cal E}_{1,4}) + \frac{1}{2}
( {\cal E}_{1,1}+ {\cal E}_{2,2}- {\cal E}_{4,4}- {\cal E}_{5,5}) \right) + \\ 
({\cal E}_{2,5} + {\cal E}_{1,4} ) \otimes \left(- \frac{1}{4}
({\cal E}_{4,5}+ {\cal E}_{1,2}) - \frac{1}{4} 
{\cal E}_{1,5} \right) + \\ ( 
{\cal E}_{1,1}+ {\cal E}_{2,2}- {\cal E}_{4,4}- {\cal E}_{5,5} ) \otimes
\left( \frac{1}{2} ( {\cal E}_{1,4}- {\cal E}_{2,5}) \right) + \\ ( 
{\cal E}_{1,1}- {\cal E}_{2,2}+ {\cal E}_{4,4}- {\cal E}_{5,5} ) \otimes
\left( \frac{1}{2} ( {\cal E}_{1,2}- {\cal E}_{4,5}) \right) + \\ ( 
{\cal E}_{4,4}- {\cal E}_{5,5}+ {\cal E}_{2,4} ) \otimes \left( \frac{1}{2} 
{\cal E}_{1,5} \right) + \\ (
{\cal E}_{1,4} ) \otimes \left( - \frac{1}{4} {\cal E}_{1,5}- {\cal E}_{2,5}
\right) + \\ ({\cal E}_{1,5} ) \otimes \left( - \frac{1}{2}
( {\cal E}_{1,1}- {\cal E}_{2,2}) + \frac{1}{2} {\cal E}_{2,4}+ 
\frac{1}{2} {\cal E}_{2,5} + \frac{1}{4} {\cal E}_{1,4}+ \frac{1}{4} 
{\cal E}_{1,5} - \frac{3}{8} {\cal E}_{4,5}+ 
\frac{5}{8} {\cal E}_{1,2} \right) . \\  
\end{array}
$$

\end{document}